\documentclass[11pt]{article}
\usepackage[utf8]{inputenc}
\usepackage{microtype}
\usepackage[a4paper,margin=1in]{geometry}
\usepackage{amsmath,amssymb,amsthm,mathtools}
\usepackage{authblk}
\usepackage{bm}
\usepackage[authoryear,round]{natbib}
\usepackage[hidelinks]{hyperref}

\usepackage[T1]{fontenc}
\usepackage{lmodern}

\hypersetup{pdftitle={Sharp Local Entropy Bounds in the Wasserstein Space},pdfauthor={Chang Jun Im and Jeong Min Jeon}}
\newtheorem{theorem}{Theorem}[section]
\newtheorem{corollary}[theorem]{Corollary}
\theoremstyle{remark}
\newtheorem{remark}[theorem]{Remark}
\theoremstyle{plain}
\numberwithin{equation}{section}
\newcommand{\W}{\mathcal W}
\newcommand{\Wc}{\mathcal{W}_{\rm cont}}
\newcommand{\Q}{\mathcal Q}
\newcommand{\M}{\mathcal M}
\newcommand{\bM}{\mathbb M}

\newcommand{\eps}{\varepsilon}
\newcommand{\dd}{\,\mathrm d}

\newcommand{\E}{\mathbb E}
\newcommand{\R}{\mathbb{R}}
\DeclareMathOperator*{\argmin}{\arg\min}
\title{Sharp Local Entropy Bounds in the Wasserstein Space}
\author[1]{Chang Jun Im}
\author[2]{Jeong Min Jeon\thanks{Corresponding author. Email: \href{mailto:jm.jeon@snu.ac.kr}{jm.jeon@snu.ac.kr}.}}
\affil[1]{The Institute for Data Innovation in Science, Seoul National University, Seoul, Korea}
\affil[2]{Department of Statistics and School of Transdisciplinary Innovations, Seoul National University, Seoul, Korea}
\date{}

\begin{document}
\maketitle

\begin{abstract}
Fr\'echet regression extends regression analysis to metric space-valued
responses such as probability distributions. Its asymptotic theory often
relies on local entropy conditions, which control the complexity of small
balls in the response space. A commonly used condition requires the
associated local entropy integral to remain bounded as the balls shrink, but
this requirement can fail when the response space is infinite-dimensional.
We determine the exact rate at which this integral diverges in the
Wasserstein space of probability distributions
supported on a compact interval. Specifically, for a broad class of ball
centers, we establish matching upper and lower bounds for the local
covering entropy and show that the local entropy integral grows at the rate of
the reciprocal square root of the ball radius. This rate persists when the
response space is restricted to a large class of continuous distributions.
The same rate holds for uniform local entropy integrals, with the supremum taken
over ball centers. These results quantify the local complexity of the
Wasserstein space and provide sharp benchmarks for entropy assumptions in
statistical problems with distribution-valued responses.
\end{abstract}



\section{Introduction}\label{sec:introduction}

Fr\'echet regression describes the relationship between a response variable \(Y\)
taking values in a general metric space \((\bM,d_\bM)\) and a predictor
\(X\) through a conditional Fr\'echet mean. The conditional Fr\'echet mean
is defined by
\[
m_\oplus(x):=\argmin_{y\in\bM}\E[d_\bM^2(Y,y)| X=x],
\]
provided that the minimizer exists and is unique. \citet{PM2019}
introduced global and local Fr\'echet regression for Euclidean predictors,
and \citet{CM2022} established uniform convergence results for local
Fr\'echet regression. Subsequent developments include functional models
for time-varying random objects \citep{DM2020}, variable selection
\citep{TWM2023}, single-index regression \citep{BM2023}, errors-in-variables
regression \citep{SH2023}, nonlinear global regression \citep{BLX2025}, partially global
regression \citep{TW2025}, deep learning methods
\citep{IZM2025}, and local constant and local linear methods for
non-Euclidean predictors \citep{Sphere2025,Torus2026}.

The asymptotic analysis of Fr\'echet regression estimators uses empirical
process arguments that require local entropy conditions on the response
space. For \(y\in\bM\) and \(r>0\), let
\(B_\bM(y,r):=\{z\in\bM:d_\bM(y,z)<r\}\). For any nonempty subset
\(S\subseteq\bM\), let \(N(r,S,d_\bM)\) denote the smallest number of open
balls of radius \(r\), with centers in \(S\), needed to cover \(S\). A commonly used condition requires the local entropy integral at a target
object \(y_\oplus\in\bM\) to remain bounded, that is,
\begin{equation}\label{eq:intro-entropy}
\int_0^1\sqrt{1+\log N\bigl(\delta\eps,B_\bM(y_\oplus,\delta),d_\bM\bigr)}\dd\eps
=O(1)\quad\text{as}\quad\delta\to0,
\end{equation}
where the target object \(y_\oplus\) is \(m_\oplus(x)\) or another object,
depending on the regression model.
Uniform versions of \eqref{eq:intro-entropy}, which take a supremum over a set of ball centers \(y_\oplus\) either inside or outside the integral, are also commonly used; see the aforementioned literature. Such local entropy integrals also appear beyond Fr\'echet regression,
for example in Fr\'echet analysis of variance \citep{DM2019} and change-point
detection \citep{DMC2020}.

In a Euclidean response space, the covering number $N(\delta\eps,B_\bM(y_\oplus,\delta),d_\bM)$ in \eqref{eq:intro-entropy}
does not depend on \(\delta\), and the condition holds. In infinite-dimensional response spaces, however, the condition may fail. A
relaxed form of \eqref{eq:intro-entropy}, used by \citet{Sphere2025}, requires only that
\begin{equation}\label{eq:intro-relaxed}
\int_0^1\sqrt{1+\log N\bigl(\delta\eps,B_\bM(m_\oplus(x),\delta),d_\bM\bigr)}\dd\eps
=O(\delta^{\alpha-1})\quad\text{as}\quad\delta\to0
\end{equation}
for some \(0<\alpha\le1\). The exponent \(\alpha\) enters the convergence
rates of their local regression estimators. In fact, under certain
conditions, one can achieve the convergence rate
\(O_p(n^{-2/(8-4\alpha+d)})\), where \(d\ge1\) is the dimension of the
predictor. Thus, a larger \(\alpha\) gives a faster convergence rate.
Similar pattern applies to global regression estimators \citep{PM2026}.

An important example of a response space in this literature is the Wasserstein space of probability distributions supported on a compact interval, equipped with the quadratic Wasserstein distance.
Such distribution-valued responses are common, since in many applications each observational unit
is naturally summarized by a probability distribution, as with
age-at-death distributions in demography, dose-response distributions, distributions of human physical activity levels or distributions extracted from brain imaging data \citep{PM2016,JPV2021,JLMP2022,CLM2023}.

The question for distribution-valued responses is therefore how large
\(\alpha\) can be taken in \eqref{eq:intro-relaxed}. For the Wasserstein
space, \citet{Sphere2025} proved that \eqref{eq:intro-relaxed} holds for arbitrary \(0<\alpha<1/2\). 
\citet{PM2026}, the correction note of \citet{PM2019}, further extended this result by showing that \eqref{eq:intro-relaxed} holds with \(\alpha=1/2\). \citet{PM2026} also verified the stronger condition \eqref{eq:intro-entropy} on a
finite-dimensional subset of the Wasserstein space and noted that it may
fail on the space itself. \citet{SPL2026} showed that \eqref{eq:intro-entropy} does fail when
the target object is a strictly increasing absolutely continuous
probability distribution whose density is bounded away from zero and infinity. They
obtained a lower bound of order \(\sqrt{\log(1/\delta)}\) for the local
entropy integral. The integral was thus known to diverge at a rate between
\(\sqrt{\log(1/\delta)}\) and \(\delta^{-1/2}\),
and the exact rate remained open.

In this paper, we determine the exact rate of divergence of the local
entropy integral not only in the Wasserstein space and but also in a large class of
continuous distributions. Specifically, we prove that the integral
diverges at the exact rate \(\delta^{-1/2}\) by establishing matching
upper and lower bounds for a broad class of ball centers. We also study the uniform versions of the local entropy
integral, with the supremum over ball centers taken either inside or
outside the integral. Under a mild condition on the set of ball centers,
both versions diverge at the same exact rate \(\delta^{-1/2}\). These
results follow from the exact order of the logarithm of the covering
number, \(\log N(\delta\eps,B_\bM(y_\oplus,\delta),d_\bM)\), which we also
establish. This order may be of independent interest, for example in the
analysis of posterior contraction \citep{Nguyen2013}, in the study of
empirical barycenters \citep{ACLP2020} and in binary regression with metric space-valued
predictors \citep{LinLin2025}, where covering numbers for balls in metric spaces are employed.

The remainder of this paper is organized as follows. Section~\ref{sec:full} establishes the local entropy bounds in the Wasserstein space. Section~\ref{sec:continuous} provides similar results for a class of continuous distributions, and Section~\ref{sec:uniform} treats the uniform versions.

\section{Pointwise local entropy bounds on Wasserstein space}\label{sec:full}

For $-\infty<a<b<\infty$, let \(\W([a,b])\) be the set of distribution functions of probability measures supported on \([a,b]\). Without loss of generality, we focus on \(\W:=\W([0,1])\). The quantile function of \(F\in\W\) is given by
\[
q_F(t):=\inf\{u\in\R:F(u)\ge t\}\in[0,1],\quad 0<t<1.
\]
For \(F,G\in\W\), the quadratic Wasserstein distance $d_2$ is given by
\begin{equation}\label{eq:quantile}
d_2(F,G):=\left(\int_0^1(q_F(t)-q_G(t))^2\dd t\right)^{1/2}.
\end{equation}
Then, $(\W,d_2)$ is a totally bounded metric space \citep{PZ}. For \(F\in\M\subseteq\W\) and \(\delta>0\), write
\begin{align*}
B_{\M}(F,\delta)&:=\{G\in\M:d_2(F,G)<\delta\},\\
I_{\M}(F,\delta)&:=\int_0^1\sqrt{1+\log N(\delta\eps,B_{\M}(F,\delta),d_2)}\dd\eps.
\end{align*}
Note that every probability measure $\mu$ supported on $[0,1]$ admits the Lebesgue decomposition
\[
\mu(A)=\int_Ag_\mu(x)\dd x+\mu^s(A),\quad A\subseteq[0,1]\text{ is Borel},
\]
where $g_\mu$ is a nonnegative measurable function on $[0,1]$, and \(\mu^s\) is a measure
singular with respect to Lebesgue measure.
Define 
\[
\mathrm{ACC}([0,1]):=\left\{\mu:\mu\text{ is a probability measure on }[0,1]\text{ such that }\int_0^1g_\mu(x)\dd x>0\right\}.
\]
In other words, $\mathrm{ACC}([0,1])$ is the set of probability measures on $[0,1]$ that have a nonzero absolutely continuous component. In particular, it includes continuous probability measures and mixture of continuous and discrete probability measures but excludes discrete probability measures.
For \(F\in\W\), let \(\mu_F\) be its probability measure. We now state the exact order of $\log N(\delta\eps,B_{\W}(F,\delta),d_2)$ and the exact rate of divergence of $I_{\W}(F,\delta)$.

\begin{samepage}
\begin{theorem}\label{thm:full}
There exists a universal constant \(C>0\) such that
\begin{equation}\label{eq:full-upper}
\log N(\delta\eps,B_{\W}(F,\delta),d_2)\le\frac{C}{\delta\eps},
\quad F\in\W,\quad 0<\delta,\eps\le1.
\end{equation}
If \(\mu_F\in\mathrm{ACC}([0,1])\), then there also exist constants \(c_F>0\) and \(0<\eps_F\le1\), depending only on \(F\), such that
\begin{equation}\label{eq:full-lower}
\log N(\delta\eps,B_{\W}(F,\delta),d_2)\ge\frac{c_F}{\delta\eps},
\quad 0<\delta\le1,\quad 0<\eps\le\eps_F.
\end{equation}
Consequently, at every such \(F\in\W\),
\[
I_{\W}(F,\delta)\asymp\delta^{-1/2}\quad\text{as}\quad\delta\to0.
\]
\end{theorem}
\end{samepage}

\begin{proof}
For the upper bound, let \(\Q\) be the set of nondecreasing functions from \((0,1)\) to \([0,1]\), regarded as elements of \(L^2(0,1)\) equipped with the norm $\|\cdot\|_2$. Theorem~1.1 of \citet{GW}, applied to functions of one variable, gives
\[
\log N_{[]}\left(\frac{\delta\eps}{2},\Q,\|\cdot\|_2\right)\le\frac{2C_0}{\delta\eps},
\quad 0<\delta,\eps\le1,
\]
where \(N_{[]}\) is the bracketing number and \(C_0>0\) is a universal constant. Since
\[
\log N(\delta\eps,S,\|\cdot\|_2)
\le\log N_{[]}\left(\frac{\delta\eps}{2},\Q,\|\cdot\|_2\right)
\]
for every nonempty subset \(S\subseteq\Q\), we obtain
\begin{equation}\label{eq:subset-upper}
\log N(\delta\eps,S,\|\cdot\|_2)
\le\frac{2C_0}{\delta\eps}.
\end{equation}
Taking \(S=\{q_G:G\in B_{\W}(F,\delta)\}\) and using \eqref{eq:quantile} proves \eqref{eq:full-upper} with \(C=2C_0\). Moreover,
\begin{equation}\label{eq:integral-upper}
I_{\W}(F,\delta)
\le1+\sqrt{\frac{C}{\delta}}\int_0^1\eps^{-1/2}\dd\eps=1+2\sqrt C\,\delta^{-1/2}
\le(1+2\sqrt C)\delta^{-1/2},\quad 0<\delta\le1.
\end{equation}

For \(m\in\mathbb N\), partition \([0,1]\) into
\[
I_j:=[(j-1)/m,j/m),\quad 1\le j<m,
\quad I_m:=[(m-1)/m,1].
\]
Set \(p_j:=\mu_F(I_j)\), \(P_j:=\sum_{\ell=1}^jp_\ell\) and \(P_0:=0\). For \(\theta=(\theta_1,\ldots,\theta_m)\in\Gamma_m:=\{0,1\}^m\), define
\begin{equation}\label{eq:polynomial-cdf}
F_\theta(x):=P_{j-1}+p_j(mx-j+1)^{1+\theta_j},
\quad x \in I_j,
\end{equation}
and extend it by \(0\) for \(x<0\) and \(1\) for \(x>1\). This function is continuous and piecewise polynomial on a finite partition. Hence, \(F_\theta\) is absolutely continuous on $[0,1]$.
If every \(p_j\) is positive, then \(F_\theta\) is strictly increasing on \([0,1]\). The quantiles of \(F\) and \(F_\theta\) satisfy
\[
q_F(t)\in I_j,\quad
q_{F_\theta}(t)\in
\left(\frac{j-1}{m},\frac jm\right),
\quad P_{j-1}<t<P_j
\]
whenever \(p_j>0\). Consequently,
\[
|q_F(t)-q_{F_\theta}(t)|<\frac1m,
\quad P_{j-1}<t<P_j
\]
whenever \(p_j>0\). Since intervals with \(p_j=0\) have \(P_j=P_{j-1}\), we obtain
\begin{equation}\label{eq:bin-radius}
d_2^2(F_\theta,F)=\int_0^1|q_{F_\theta}(t)-q_F(t)|^2\dd t=\sum_{j:p_j>0}\int_{P_{j-1}}^{P_j}
  |q_{F_\theta}(t)-q_F(t)|^2\dd t<\frac1{m^2}\sum_{j=1}^mp_j=\frac1{m^2}.
\end{equation}

For $\theta,\varphi\in\Gamma_m$ and \(j\) with \(p_j>0\), \eqref{eq:polynomial-cdf} gives
\[
q_{F_\theta}(t)
=\frac{j-1+u_{j,t}^{1/(1+\theta_j)}}{m},
\quad
q_{F_{\varphi}}(t)
=\frac{j-1+u_{j,t}^{1/(1+\varphi_j)}}{m},\quad P_{j-1}<t<P_j,
\]
where $u_{j,t}=(t-P_{j-1})/p_j$.
Since \(\theta_j,\varphi_j\in\{0,1\}\), we obtain
\[
|q_{F_\theta}(t)-q_{F_{\varphi}}(t)|^2
=\frac{|\theta_j-\varphi_j|}{m^2}(u_{j,t}-\sqrt{u_{j,t}})^2,\quad P_{j-1}<t<P_j,
\]
which implies that
\begin{equation}\label{eq:distance-identity}
\begin{aligned}
d_2^2(F_\theta,F_{\varphi})&=\int_0^1
  |q_{F_\theta}(t)-q_{F_{\varphi}}(t)|^2\dd t=\sum_{j:p_j>0}\int_{P_{j-1}}^{P_j}
  |q_{F_\theta}(t)-q_{F_{\varphi}}(t)|^2\dd t\\
  &=\frac1{30m^2}\sum_{j=1}^m
  p_j|\theta_j-\varphi_j|.
\end{aligned}
\end{equation}

To select distributions from $\{F_\theta\}$ that are sufficiently far apart in Wasserstein distance, fix a measurable function \(h\) on \([0,1]\) and a constant \(K\ge1\) such that
\begin{equation*}
\int_A h(x)\dd x\le\mu_F(A),\quad A\subseteq[0,1]\text{ is Borel},\quad
0\le h\le K\text{ a.e.},\quad a:=\int_0^1h(x)\dd x>0.
\end{equation*}
Such $h$ and $K$ exist since \(\mu_F\in\mathrm{ACC}([0,1])\). For example, one may take \(h=g_{\mu_F}\wedge1\) and \(K=1\). Set
\[
w_j:=\int_{I_j}h(x)\dd x.
\]
Then,
\[
w_j\le p_j,\quad w_j\le K/m,\quad \sum_{j=1}^mw_j=a.
\]
For any set $D$, let $|D|$ denote its cardinality. For \(\theta\in\Gamma_m\),
\begin{align*}
&\left|\left\{\varphi\in\Gamma_m:
  \sum_{j=1}^m w_j|\theta_j-\varphi_j|\le a/4
  \right\}\right|\\
&\quad=
  \sum_{\varphi\in\Gamma_m}
  \mathbf1_{\{\sum_{j=1}^m
  w_j|\theta_j-\varphi_j|\le a/4\}}\\
&\quad\le\sum_{\varphi\in\Gamma_m}\exp\left[
\frac mK\left(
\frac a4-\sum_{j=1}^m w_j|\theta_j-\varphi_j|
\right)
\right]\\
&\quad=
  e^{am/(4K)}
  \sum_{\varphi\in\Gamma_m}
  \exp\left(
    -\frac mK\sum_{j=1}^m
    w_j|\theta_j-\varphi_j|
  \right)\\
&\quad=
  e^{am/(4K)}
  \prod_{j=1}^m
  \left(
    \sum_{\varphi_j=0}^1
    e^{-mw_j|\theta_j-\varphi_j|/K}
  \right)\\
&\quad=
  e^{am/(4K)}
  \prod_{j=1}^m(1+e^{-mw_j/K})\\
&\quad\le
  2^m\exp\left(
    -\frac{am}{4K}
    +\frac{m^2}{8K^2}\sum_{j=1}^m w_j^2
  \right)\\
&\quad\le 2^m e^{-am/(8K)},
\end{align*}
where the second inequality follows from
\(1+e^{-x}\le2e^{-x/2+x^2/8}\), and the last inequality follows from
\(\sum_jw_j^2\le(K/m)\sum_jw_j=Ka/m\). Choose maximal \(\Theta_m\subseteq\Gamma_m\) under inclusion
subject to
\[
\sum_{j=1}^m w_j|\theta_j-\varphi_j|>a/4,
\quad
\theta\ne\varphi,\quad\theta,\varphi\in\Theta_m.
\]
By maximality, every vector in \(\Gamma_m\) belongs to $\{\varphi\in\Gamma_m:
  \sum_{j=1}^m w_j|\theta_j-\varphi_j|\le a/4\}$ for some \(\theta\in\Theta_m\).
Consequently,
\[
2^m\le|\Theta_m|\,2^m e^{-am/(8K)},
\]
which implies that
\begin{equation*}
\log|\Theta_m|\ge\frac{am}{8K}.
\end{equation*}

Set
\begin{equation*}
\eps_F:=\frac14\sqrt{\frac a{120}},
\quad
c_F:=\frac{a\eps_F}{6K}.
\end{equation*}
For \(0<\delta\le1\) and \(0<\eps\le\eps_F\), choose
\[
m:=\left\lfloor\frac{2\eps_F}{\delta\eps}\right\rfloor.
\]
Since \(2\eps_F/(\delta\eps)\ge2/\delta\ge2\), the inequality
\(\lfloor t\rfloor\ge 2t/3\) for \(t\ge2\) gives
\begin{equation*}
\frac{4\eps_F}{3\delta\eps}
\le m\le\frac{2\eps_F}{\delta\eps},
\quad
\frac1m\le\frac{3\delta\eps}{4\eps_F}
\le\frac{3\delta}{4}<\delta.
\end{equation*}
This together with \eqref{eq:bin-radius} implies that \(F_\theta\in B_{\W}(F,\delta)\) for every
\(\theta\in\Theta_m\).
By \eqref{eq:distance-identity} and \(p_j\ge w_j\), distinct
\(\theta,\varphi\in\Theta_m\) satisfy
\begin{equation*}
d_2(F_\theta,F_\varphi)
>\frac{4\eps_F}{m}
\ge2\delta\eps.
\end{equation*}
Thus, any ball of radius \(\delta\eps\) contains at most one
\(F_\theta\) with \(\theta\in\Theta_m\). Consequently,
\begin{equation}\label{eq:explicit-lower}
\log N(\delta\eps,B_{\W}(F,\delta),d_2)
\ge\log|\Theta_m|
\ge\frac{am}{8K}
\ge\frac{a\eps_F}{6K\delta\eps}
=\frac{c_F}{\delta\eps}.
\end{equation}
This proves \eqref{eq:full-lower}. 

Since
\(0<a\le1\) and \(K\ge1\), we obtain
\[
1+\log N(\delta\eps,B_{\W}(F,\delta),d_2)
\ge1+\frac{am}{8K}
\ge\frac{a(m+1)}{8K}
>\frac{a\eps_F}{4K\delta\eps},
\]
where the last inequality follows from
\(m+1>2\eps_F/(\delta\eps)\).
Therefore,
\begin{equation}\label{eq:integral-lower}
I_{\W}(F,\delta)
\ge\int_0^{\eps_F}
  \sqrt{\frac{a\eps_F}{4K\delta\eps}}\,\dd\eps
=\eps_F\sqrt{\frac aK}\,\delta^{-1/2}
=\frac{a}{4\sqrt{120K}}\,\delta^{-1/2}.
\end{equation}
This together with \eqref{eq:integral-upper} proves
\(I_{\W}(F,\delta)\asymp\delta^{-1/2}\)
as \(\delta\to0\).
\end{proof}

\section{Class of strictly increasing absolutely continuous distributions}\label{sec:continuous}

In applications, distribution-valued responses often take values in the following class of continuous distributions:
\[
\Wc:=\bigl\{F\in\W:F(0)=0,\ 
F\text{ is strictly increasing and absolutely continuous on }[0,1]\bigr\}.
\]
Note that every $F\in\Wc$ satisfies that \(\mu_F\in\mathrm{ACC}([0,1])\).
One might expect that this class yields a slower divergence of
the local entropy integral as \(\delta\to0\). The following theorem
shows that the rate of divergence remains \(\delta^{-1/2}\), even when
both the ball centers and the local balls are restricted to \(\Wc\).

\begin{theorem}\label{thm:continuous}
There exists a universal constant \(C>0\) such that
\[
\log N(\delta\eps,B_{\Wc}(F,\delta),d_2)\le\frac{C}{\delta\eps},
\quad F\in\Wc,\quad 0<\delta,\eps\le1.
\]
There also exist constants \(c_F>0\) and \(0<\eps_F\le1\), depending only on \(F\in\Wc\), such that
\[
\log N(\delta\eps,B_{\Wc}(F,\delta),d_2)\ge\frac{c_F}{\delta\eps},
\quad 0<\delta\le1,\quad 0<\eps\le\eps_F.
\]
Consequently, at every $F\in\Wc$,
\[
I_{\Wc}(F,\delta)\asymp\delta^{-1/2}\quad\text{as}\quad\delta\to0.
\]
\end{theorem}

\begin{proof}
Recall the notation in the
proof of Theorem~\ref{thm:full}. Applying \eqref{eq:subset-upper} to
\[S=\{q_G:G\in B_{\Wc}(F,\delta)\}\] and using \eqref{eq:quantile}, we obtain the desired upper covering bound for \(\log N(\delta\eps,B_{\Wc}(F,\delta),d_2)\). By arguing as in
the proof of \eqref{eq:integral-upper}, we obtain
\[
I_{\Wc}(F,\delta)
\le(1+2\sqrt C)\delta^{-1/2},\quad 0<\delta\le1.
\]

For the lower bound, strict increase of \(F\) on $[0,1]$ ensures that every \(p_j\) is positive. Each \(F_\theta\) in
\eqref{eq:polynomial-cdf} is therefore strictly increasing on
\([0,1]\). Since \(F_\theta\) is absolutely continuous on
\([0,1]\) and satisfies \(F_\theta(0)=0\), it belongs to \(\Wc\).
By arguing as in the proof of Theorem~\ref{thm:full}, we obtain
the lower bounds \eqref{eq:explicit-lower} and
\eqref{eq:integral-lower} with \(\W\) replaced by \(\Wc\).
\end{proof}

\section{Uniform local entropy bounds}\label{sec:uniform}

For \(\M\in\{\W,\Wc\}\), a nonempty subset
\(\mathcal A\subseteq\M\), and \(\delta>0\), define
\begin{align*}
I_{\M}^{\mathrm{out}}(\mathcal A,\delta)
&:=\sup_{F\in\mathcal A}I_{\M}(F,\delta),\\
I_{\M}^{\mathrm{in}}(\mathcal A,\delta)
&:=\int_0^1\sup_{F\in\mathcal A}
\sqrt{1+\log N(\delta\eps,B_{\M}(F,\delta),d_2)}\dd\eps.
\end{align*}
To accommodate \(\W\)-valued responses, \citet{PM2019} assumed
that \(I_{\W}^{\mathrm{in}}(\mathcal A,\delta)=O(1)\) as
\(\delta\to0\), with
\(\mathcal A=\{m_{\rm glob}(x):x\in\mathbb R^p,\ \|x\|<r\}\)
for some \(r>0\), where $m_{\rm glob}(x)$ is the global Fr\'echet regression function.
Similarly, \citet{Torus2026} assumed that
\(I_{\W}^{\mathrm{out}}(\mathcal A,\delta)
=O(\delta^{\alpha-1})\) as \(\delta\to0\) for some
\(0<\alpha\le1\), with
\(\mathcal A=\{F\in\W:d_2(m_{\oplus}(x),F)<r\}\)
for some \(r>0\). Therefore, it is important to determine the rates of divergence of these
uniform local entropy integrals. The following theorem
shows that these integrals still diverge at the rate \(\delta^{-1/2}\).

\begin{theorem}\label{thm:subset-suprema}
There exists a universal constant \(C_I>0\) such that, for every
\(\M\in\{\W,\Wc\}\) and every nonempty subset
\(\mathcal A\subseteq\M\),
\begin{equation}\label{eq:upper bound}
I_{\M}^{\mathrm{out}}(\mathcal A,\delta)
\le I_{\M}^{\mathrm{in}}(\mathcal A,\delta)
\le C_I\,\delta^{-1/2},
\quad 0<\delta\le1.
\end{equation}
If \(\mathcal A\) contains a distribution \(H\) with
\(\mu_H\in\mathrm{ACC}([0,1])\), then there exists a constant
\(c_{\mathcal A}>0\), depending only on $\mathcal A$, such that
\begin{equation}\label{eq:lower bound}
I_{\M}^{\mathrm{out}}(\mathcal A,\delta)
\ge c_{\mathcal A}\,\delta^{-1/2},
\quad 0<\delta\le1.
\end{equation}
In particular, if \(\mathcal A\) contains an absolutely continuous
distribution with a density bounded by \(K\ge1\) almost everywhere, then the lower bound in \eqref{eq:lower bound} holds with
\begin{equation}\label{eq:bounded-density-constant}
c_{\mathcal A}=\frac1{4\sqrt{120K}}.
\end{equation}
\end{theorem}

\begin{proof}
Recall the notation in the
proof of Theorem~\ref{thm:full}. Applying \eqref{eq:subset-upper} to
\[S=\{q_G:G\in B_{\M}(F,\delta)\}\] and using
\eqref{eq:quantile}, we obtain
\[
\log N(\delta\eps,B_{\M}(F,\delta),d_2)
\le\frac{C}{\delta\eps},\quad
F\in\M,\quad 0<\delta,\eps\le1,
\]
where \(C>0\) is a universal constant. Hence,
\[
I_{\M}^{\mathrm{out}}(\mathcal A,\delta)
\le I_{\M}^{\mathrm{in}}(\mathcal A,\delta)\le\int_0^1\sqrt{1+\frac{C}{\delta\eps}}\dd\eps
\le(1+2\sqrt C)\delta^{-1/2},\quad 0<\delta\le1.
\]
Thus, \(C_I:=1+2\sqrt C\) is a desired universal constant.

For the lower bound, choose \(H\in\mathcal A\) with
\(\mu_H\in\mathrm{ACC}([0,1])\). By arguing as in the proof of \eqref{eq:integral-lower}, we obtain the lower bound
\[
I_{\M}^{\mathrm{out}}(\mathcal A,\delta)
\ge I_{\M}(H,\delta)
\ge\frac{a}{4\sqrt{120K}}\,\delta^{-1/2},
\quad 0<\delta\le1.
\]

If \(H\) has a density bounded by \(K\ge1\) almost everywhere,
use that density as \(h\), so that \(a=1\). The preceding
bound then yields
\[
I_{\M}(H,\delta)
\ge\frac1{4\sqrt{120K}}\,\delta^{-1/2}.
\]
\end{proof}

Note that every nonempty \(\mathcal A\subseteq\Wc\) contains a distribution \(H\) with
\(\mu_H\in\mathrm{ACC}([0,1])\). A similar upper bound to \eqref{eq:upper bound} was obtained by \citet{PM2026} for $\M=\W$, but lower bounds were not covered. We now apply Theorem~\ref{thm:subset-suprema} to \[\mathcal A=\{F\in\M:d_2(F_\oplus,F)<r\}=:\mathcal S_\M(F_\oplus,r)\] for \(\M\in\{\W,\Wc\}\), \(F_\oplus\in\M\) and \(r>0\). This gives the exact rate of divergence of the integral in the assumption of \citet{Torus2026} when $\M=\W$. In particular, \eqref{eq:bounded-density-constant} allows the constant in the lower bound to be chosen uniformly in \(F_\oplus\).

\begin{corollary}\label{cor:neighborhood}
There exist a constant \(c_r>0\),
depending only on \(r\), and a universal constant \(C_I>0\) such that, for every
\(\M\in\{\W,\Wc\}\),
\[
c_r\delta^{-1/2}\le I_{\M}^{\mathrm{out}}(\mathcal S_\M(F_\oplus,r),\delta)\le I_{\M}^{\mathrm{in}}(\mathcal S_\M(F_\oplus,r),\delta)
\le C_I\delta^{-1/2},\quad F_\oplus\in\M,\quad 0<\delta\le1.
\]
\end{corollary}

\begin{proof}
Recall the notation in the
proof of Theorem~\ref{thm:full}. Set
\[
m_r:=\lceil1/r\rceil.
\]
Apply the construction in \eqref{eq:polynomial-cdf} with
\(F=F_\oplus\) and \(m=m_r\), and let \(G:=F_{\theta}\) for $\theta=(0,\ldots,0)$.
Then, \(G\) is absolutely continuous on \([0,1]\) and satisfies \(G(0)=0\).
If \(\M=\Wc\), then strict increase of \(F_\oplus\) on $[0,1]$ ensures that
every \(p_j\) is positive. In this case, \(G\) is strictly increasing
on \([0,1]\). Consequently, \(G\in\M\) for either choice of \(\M\).
Moreover, \eqref{eq:bin-radius} gives
\[
d_2(G,F_\oplus)<\frac{1}{m_r}\le r.
\]
Hence, \(G\in\mathcal S_{\M}(F_\oplus,r)\).
Since the density of $G$ is bounded by \(m_r\) almost
everywhere, Theorem~\ref{thm:subset-suprema} with \(K=m_r\)
yields the asserted bounds with
\[
c_r:=\frac1{4\sqrt{120m_r}}.
\]
\end{proof}

\begin{remark}
The same asymptotic orders hold for general \(\W([a,b])\). They also remain valid when the radius $\delta\eps$ of the covering balls or the radius $\delta$ of the local ball in $\log N(\delta\eps,B_{\M}(F,\delta),d_2)$ are rescaled by positive constants. A similar upper bound to \eqref{eq:upper bound} was obtained by \citet{PM2026} for $\M=\W$, but lower bounds are missing.
\end{remark}

\section*{Acknowledgements}

Chang Jun Im was supported by the National Research Foundation of Korea grant funded by the Korea government (MSIT) No. RS-2025-00515381. Jeong Min Jeon was supported by the National Research Foundation of Korea grant funded by the Korea government (MSIT) No. RS-2023-00211910.

\begingroup
\setlength{\bibsep}{5pt}

\endgroup
\end{document}